\documentclass[12pt]{article}

\usepackage{array}
\usepackage{amsmath,amssymb}
\usepackage{bm}

\newcommand{\dd}{\bm{d}}
\newcommand{\Ddd}{\bm{D}}

\newcommand{\ggt}{\bm{G}'}
\newcommand{\g}{\bm{G}}
\newcommand{\hh}{\bm{H}}

\renewcommand{\lll}{\bm{l}}
\newcommand{\Lll}{\bm{L}}

\newcommand{\qq}{\bm{q}}

\newcommand{\sss}{\bm{S}}
\newcommand{\ww}{\bm{W}}

\newcommand{\qqq}{\bm{Q}}

\newcommand{\ggamma}{\bm{\Gamma}}

\newcommand{\ggammat}{\bm{\Gamma}'}
\newcommand{\Llambda}{\bm{\Lambda}}

\newcommand{\Llambdain}{\Llambda^{-1}}
\newcommand{\llambda}{\bm{\lambda}}
\newcommand{\s}{\bm{\Sigma}}

\newcommand{\gtilde}{\widetilde{\bm{G}}}
\newcommand{\gtildet}{\gtilde'}
\newcommand{\iip}{\bm{I}_p}

\newcommand{\xxi}{\bm{\xi}}
\newcommand{\xxxi}{\bm{\Xi}}

\newcommand{\ddd}{\bm{d}}

\newcommand{\op}{{\cal O}(p)}

\newcommand{\tr}{\mathop{\rm tr}}

\newcommand{\diag}{\mathop{\rm diag}}
\newcommand{\proof}{\noindent{\bf Proof.}\quad}
\newcommand{\qed}{\hbox{\rule[-2pt]{3pt}{6pt}}}

\newcommand{\sumi}{\sum_{i=1}^p}

\newcommand{\dlim}{\stackrel{d}{\rightarrow}}

\newtheorem{lem}{Lemma}
\newtheorem{theorem}{Theorem}

\renewcommand{\aa}{\alpha}
\newcommand{\bb}{\beta}

\title{Inference on Eigenvalues of  Wishart Distribution Using Asymptotics with respect to the Dispersion of Population Eigenvalues}
\author{Yo Sheena\thanks{Department of Economics, Shinshu University}
  \ and Akimichi Takemura\thanks{Graduate School of Information Science and Technology, University of Tokyo}}
\date{April, 2007}

\begin{document}
\maketitle
\begin{abstract}
In this paper we derive some new and practical results on testing and interval
estimation problems for the population eigenvalues of a Wishart
matrix based on the asymptotic theory for block-wise infinite
dispersion of the population eigenvalues.  
This new type of
asymptotic theory has been developed by the present authors in
Takemura and Sheena (2005) and Sheena and Takemura (2007a,b)
and in these papers it was applied to 
point estimation problem of population covariance matrix 
in a decision theoretic framework.
In this paper we apply it to some testing and interval estimation problems.
We show that the 
approximation based on this type of asymptotics
is generally much better than the traditional large-sample asymptotics for
the problems.
\end{abstract}

\section{Introduction}
\label{sec:intro}
Let $\sss=(s_{ij})$ be distributed according to Wishart distribution
$\ww_p(n,\s)$, where $p$ is the dimension, $n$ is the degrees of
freedom and $\s$ is the covariance matrix.
Let  $\lambda_1 \ge\ldots\ge \lambda_p>0$ denote the eigenvalues of
$\s$.  In this paper we consider some testing and interval estimation
problems for the eigenvalues of $\s$.  
Our aim is to give practical solutions to the problems
based on the asymptotic theory
for block-wise infinite dispersion of the population eigenvalues.  
In view of the intractability of the finite sample exact distribution
of sample eigenvalues,
usually the large sample asymptotic approximation is used.  
There exists an extensive literature on improving the first-order large sample
approximation by an asymptotic expansion (see Siotani et al.\ (1985) for
a comprehensive treatment).
However for a moderate or small value of the sample size $n$, the large sample
asymptotic theory often gives a poor approximation.  In these cases
asymptotic expansions tend to give an even larger error.
On the other hand, we find that approximation based on the infinite
dispersion of the population eigenvalues is more robust and gives a
reasonable approximation even for a small sample size $n$.

The first problem we consider in this paper is testing the one-sided
null hypothesis on the $m$th population eigenvalue
\begin{equation}
\label{eq:test1}
H_0^{(m)}: \lambda_m\geq \lambda_m^* .
\end{equation}
For testing $H_0^{(m)}$
it is natural to consider a one-sided rejection region based on
the $m$th sample eigenvalue of $l_m$ of $\sss$. We show that
the least favorable distribution is given by 
$\lambda_m^*=\lambda_1 = \dots = \lambda_m$ and
$0=\lambda_{m+1}=\dots=\lambda_p$.   This is exactly the situation
covered by the asymptotic theory for block-wise infinite
dispersion.  Therefore it 
gives an explicit solution to the testing problem
of $H_0^{(m)}$.  

The second problem is the interval estimation for
the largest population eigenvalue $\lambda_1$ in terms of the largest
sample eigenvalue $l_1$ of $\sss$.  
We will show that confidence interval based 
block-wise infinite dispersion  gives much better coverage probability
than the conventional large sample asymptotics.

The third problem is testing the hypothesis of 
equality of the several smallest eigenvalues:
$\lambda_{m+1}=\cdots=\lambda_p$.
This problem is
important in determining the rank of the systematic part in a
multivariate variance component model.
We consider approximation to the null distribution of the likelihood ratio
criterion under the block-wise infinite dispersion of population
eigenvalues.  Again this type of asymptotics gives much better
approximation than the large sample asymptotics.

The organization of the paper is as follows. 
In Section \ref{eq:preliminaries} we set up notations for the paper
and give some preliminary results on the asymptotic theory for
block-wise infinite  dispersion of the population eigenvalues.  
In Section \ref{sec:three-problems} we study the above three problems, 
1) One-sided test for a population eigenvalue in Section 3.1; \ 
2) Interval estimation for extreme  eigenvalues in Section 3.2; \ 
3) Testing equality of the smallest eigenvalues in Section 3.3.


\section{Asymptotic Distribution of Normalized Sample Eigenvalues}
\label{eq:preliminaries}

Denote the spectral decompositions of $\s$ and $\sss$  by
\begin{align}
\label{spect_decomp_population}
\s&=\ggamma\Llambda\ggammat \\
\label{spect_decomp_sample}
\sss&=\g\Lll\ggt,
\end{align}
where $\g, \ggamma \in \op,$ the group of $p\times p$ orthogonal matrices, and 
$\Llambda=\diag(\lambda_1,\ldots,\lambda_p)$,
$\Lll=\diag(l_1, \ldots,l_p)$ are diagonal matrices with
the eigenvalues  
$\lambda_1 \ge\ldots\ge \lambda_p>0$, $l_1\ge\ldots\ge l_p>0$
of $\s$ and $\sss$, respectively. 
We use the notations $\llambda=(\lambda_1,\ldots,\lambda_p)$  and $\lll=(l_1,\ldots,l_p)$ hereafter. 
We guarantee the uniqueness (almost surely) of the decomposition (\ref{spect_decomp_sample}) by requiring
that 
\begin{equation}
\label{def_Gtilde}
\gtilde=(\widetilde{g}_{ij})=\ggammat\g
\end{equation}
has positive diagonal elements. 

In Takemura and Sheena (2005) we considered what happens to
appropriately normalized components of $\sss$ 
if the population eigenvalues become
infinitely dispersed, i.e.,
$$
(\lambda_2/\lambda_1, \lambda_3/\lambda_2, \ldots, \lambda_{p}/\lambda_{p-1})\rightarrow 0.
$$
In Sheena and Takemura (2007b) we generalized the
asymptotic result of Takemura and Sheena (2005) to the case when the
population eigenvalues are block-wise infinitely dispersed.


Let the population eigenvalues be parameterized as follows;
\begin{equation}
\label{def_tillamb}
\lambda_i=
\left\{
\begin{array}{cl}
\xi_i \aa,& \mbox{ if }i=1,\ldots, m,\\
\xi_i \bb,& \mbox{ if }i=m+1,\ldots, p,
\end{array}
\right.
\end{equation}
where $\xi_i$'s are fixed and ``asymptotic parameter'' $\aa$ and $\bb$ vary. 
When we say population eigenvalues are ``(two-)block-wise infinitely dispersed'', it means that
\begin{equation}
\bb/\aa \rightarrow 0.
\end{equation}
The above notation is used as a general notation including specific convergences (divergences) such as $(\aa, \bb)\rightarrow (\infty, 1)$, 
$(\aa, \bb)\rightarrow (1, 0)$ and so on. 
More precisely, the operation $\lim_{\bb/\aa \rightarrow 0}f(\aa,\bb)$
means  $\lim_{i\rightarrow \infty}f(\aa_i, \bb_i)$ with any specific
sequences $\aa_i, \bb_i,\ i=1, 2, \ldots$, 
such that $\bb_i/\aa_i\rightarrow 0$ as $i \rightarrow \infty.$

As Sheena and Takemura (2007b) indicates,  appropriate normalization for the sample eigenvalues is given by
\begin{equation}
d_i=
\left\{
\begin{array}{cl}
\label{def_till}
l_i /\aa,& \mbox{ if }i=1,\ldots, m,\\
l_i/ \bb,& \mbox{ if }i=m+1,\ldots, p,
\end{array}
\right.
\end{equation}
while (\ref{def_Gtilde}) also serves as an appropriate normalization for sample eigenvectors. 

For the normalized population or sample eigenvalues, we use the following notations; 
\begin{align*}
\xxi_1&=(\xi_1,\ldots,\xi_m), &\xxi_2&=(\xi_{m+1},\ldots,\xi_p),\\
\dd_1&=(d_1,\ldots,d_m), &\dd_2&=(d_{m+1},\ldots,d_p),\\
\xxxi_1&=\diag(\xi_1,\ldots,\xi_m),&\xxxi_2&=\diag(\xi_{m+1},\ldots,\xi_p),\\
\Ddd_1&=\diag(d_1,\ldots,d_m), &\Ddd_2&=\diag(d_{m+1},\ldots,d_p).
\end{align*}
Now we state the basic theorem on the asymptotic distributions of $\dd_1, \dd_2$. 
\begin{theorem}
\label{basic_theorem}
Suppose that we have two independent Wishart distributions
$$
\label{asympt_dist_lemma1}
\widetilde{\ww}_{11} \sim \ww_m(n,\xxxi_1),\nonumber\qquad
\widetilde{\ww}_{22} \sim \ww_{p-m}(n-m,\xxxi_2)\\ 
$$
and that their spectral decompositions are given by
\begin{align*}
\widetilde{\ww}_{11}&=\gtilde_{11}\widetilde{\Ddd}_1\gtildet_{11}, & 
\widetilde{\Ddd}_1&=\diag(\widetilde{d}_1,\ldots,\widetilde{d}_m),& 
\widetilde{\ddd}_1&=(\widetilde{d}_1,\ldots,\widetilde{d}_m),\\
\widetilde{\ww}_{22}&=\gtilde_{22}\widetilde{\Ddd}_2\gtildet_{22}, &
\widetilde{\Ddd}_2&=\diag(\widetilde{d}_{m+1},\ldots,\widetilde{d}_p),&
\widetilde{\ddd}_2&=(\widetilde{d}_{m+1},\ldots,\widetilde{d}_p),
\end{align*}
where  $\gtilde_{11} \in {\cal O}(m),\ \gtilde_{22}\in {\cal O}(p-m),\  
\widetilde{d}_1\geq \cdots \geq \widetilde{d}_m,\ \widetilde{d}_{m+1}\geq \cdots \geq \widetilde{d}_p.$
Then as $\bb/\aa\rightarrow 0$, 
$$
\dd_i\dlim \widetilde{\dd}_i,\quad i=1,2.
$$
\end{theorem}
\proof
Using Lemma 1 of Sheena and Takemura (2007b), we prove the convergence of the moment generating function. Let
$$
x(\g,\lll,\llambda,\alpha,\beta)=\exp\left(\alpha^{-1}\sum_{i=1}^ml_i\theta_i+\beta^{-1}\sum_{i=m+1}^pl_i\theta_i\right)
=\exp\left(\sumi d_i\theta_i\right),
$$
where $|\theta_i|<3^{-1}\min_j \xi_j^{-1}$, $\forall i$.
Notice that (19) in Lemma 1 of Sheena and Takemura (2007b) is satisfied since
\begin{align*}
x(\ggamma\g,\lll,\llambda,\alpha,\beta)
&\leq\exp\left(\alpha^{-1}\sum_{i=1}^m l_i|\theta_i|+\beta^{-1}\sum_{i=m+1}^p l_i|\theta_i|\right)\\
&\leq\exp\left(3^{-1}\alpha^{-1}\sum_{i=1}^m l_i\xi_i^{-1}+3^{-1}\beta^{-1}\sum_{i=m+1}^p l_i\xi_i^{-1}\right)\\
&=\exp\left(3^{-1}\sumi l_i \lambda_i^{-1}\right)\\
&\leq \exp\left(\tr 3^{-1}\g\Lll\ggt\Llambdain\right),
\quad \forall \g \in \op,\ \forall \lll \in \{\lll|l_1\geq \cdots \geq l_p\geq 0\}.
\end{align*}
For the last inequality, see e.g. Marshall and Olkin (1979) Ch.20.A.1. 
Since $x(\dd,\qq,\xxi,\alpha,\beta;\ggamma,\hh^{(\tau)})=\exp\left(\sumi d_i\theta_i\right)$, trivially we have
$$
\bar{x}_{\ggamma}(\hh^{(\tau)}\g(\qq_{11},\qq_{22},\bm{0}),\dd,\qqq_{21},\xxi)=\exp\left(\sumi d_i\theta_i\right).
$$
Therefore we have
\begin{align*}
\lim_{\beta/\alpha \rightarrow 0} E\left[\exp\left(\sumi d_i\theta_i \right)\right]
&=E\left[\exp\left(\sum_{i=1}^m \widetilde{d}_i(\widetilde{\ww}_{11})\theta_i+\sum_{i=m+1}^{p} \widetilde{d}_i(\widetilde{\ww}_{22})\theta_i\right)\right].
\end{align*}
\hfill \qed

%
%
%
%
%
%
%
\section{Inference on Population Eigenvalues}
\label{sec:three-problems}
The asymptotic result in the previous section has possibly various applications for inference on the population eigenvalues. We give three inference problems as interesting applications.
%
%
%
%
%
%
\subsection{One-sided Test for Population Eigenvalue}
\label{subsec3-1}
Consider the null hypothesis on the $m$th ($m=1,\ldots,p$) population eigenvalue 
$$
H_0^{(m)}: \lambda_m\geq \lambda_m^*
$$
against the alternative $H_1^{(m)}: \lambda_m < \lambda_m^*$. Need for testing $H_0$ arises in some practical cases, for example:
\begin{itemize}
\item In principal component analysis, 
$\lambda^*(=\lambda_1^*=\cdots =\lambda_p^*)$ may be 
a cut-off value and a test for $H_0^{(m)}$ is  repeatedly carried out starting from $m=1$ until $H_0^{(m)}$ is rejected. This is one of the methods for  deciding the dimension of the principal components.
\item Let $x_i (i=1,\ldots,p)$ be the return of the $i$th asset in finance and $\bm{x}=(x_1,\ldots,x_p)$ is distributed as the $p$-dimensional normal distribution $N_p(\bm{0},\s)$. $H_1^{(1)}$ is equivalent to the assertion $\bm{a}'\s \bm{a}< \lambda_1^*,\ \forall \bm{a}=(a_1,\ldots,a_p)$ such that $\|\bm{a}\|=1$. If $H_0^{(1)}$ is rejected, then it means that the group of assets $\bm{x}$ is stable in view of volatility since any portfolio among the group is never beyond 
$\lambda_1^*$ in its variance.
\end{itemize}
A natural rejection region in testing $H_0^{(m)}$ is given by $l_m\leq l_m^*(\gamma)$ for a given significance level $\gamma$. The following lemma and Theorem \ref{basic_theorem} give the critical point $l_m^*(\gamma)$. 
\begin{lem}
For any positive $c$
$$
\sup_{H_0^{(m)}}P_{\Llambda}(l_m \leq c)=\lim_{\beta\rightarrow 0}P_{\bar{\Llambda}}(l_m \leq c),
$$
where $\bar{\Llambda}=\diag(\bar{\lambda}_1,\ldots,\bar{\lambda}_p)$,  
$
\bar{\lambda}_1=\cdots=\bar{\lambda}_m=\lambda_m^*,\  \bar{\lambda}_{m+1}=\cdots=\bar{\lambda}_p=\beta.
$
\end{lem}
\proof
According to Theorem 1 of Anderson and Das Gupta (1964),  $P_{\Llambda}(l_m \leq c)$ is a monotonically decreasing function with respect to each $\lambda_i,\ (i=1,\ldots,p)$, hence 
$$
P_{\Llambda}(l_m \leq c)\leq P_{\bar{\Llambda}}(l_m \leq c),
$$
where $\beta=\lambda_p$.
Furthermore $P_{\bar{\Llambda}}(l_m \leq c)$ is monotonically increasing as $\beta$ goes to zero.
\hfill \qed

\bigskip
Because of the result of Theorem \ref{basic_theorem} with $\alpha=1$, $\xi_i=\lambda_m^*,\ (i=1,\ldots,m)$ and $\xi_i=1,\ (i=m+1,\ldots,p)$,  
$$
\lim_{\beta\rightarrow 0}P_{\bar{\Llambda}}(l_m \leq c)=P(\widetilde{l}_m\leq c),
$$
where $\widetilde{l}_m$ is distributed as the smallest eigenvalues of
$\ww_m(n, \lambda_m^*\bm{I}_m)$.  Therefore we have the following result.
\begin{theorem}
For testing hypothesis $H_0^{(m)}$ against $H_1^{(m)}$, a test with significance level $\gamma$ is given with the rejection region
$$
l_m\leq l_m^*(\gamma),
$$
where $l_m^*(\gamma)$ is the lower 100$\gamma$\% point of the smallest eigenvalue of $\ww_m(n, \lambda_m^*\bm{I}_m).$
\end{theorem}
For analytic calculation of  $l_m^*(\gamma)$, see Thompson (1962),
Hanumara and Thompson (1968). In the case $m=1$, which is practically the most important, it is given by $\lambda_1^*\chi^2_n(\gamma),$ where $\chi^2_n(\gamma)$ is the lower $100\gamma$\% point of the $\chi^2$ distribution with the degree of freedom $n.$

%
%
%
%
%
%
\subsection{Interval Estimation of Extreme Eigenvalues}
\label{subsec3-2}
In this subsection we present a new way of constructing  a confidence interval for the extreme population eigenvalues.
Let $\lambda_1\leq f_1(\lll)$ be a one-sided estimated interval with confidence level $\gamma$. For example, in the second case in Section \ref{subsec3-1}, the maximum volatility in all possible portfolio among the assets $\bm{x}$ is estimated to be less than or equal to $f_1(\lll)$. 

However if we use the exact finite distribution theory, it is not easy to find an appropriate $f_1(\lll)$ under a given $\gamma$ even if we only consider an interval of the simplest form 
$\lambda_1 \leq c_1l_1$ with some constant $c_1$. (Note that
$l_i/\lambda_i,\ (i=1,\ldots,p)$ is bounded in probability. See Lemma
1 of Takemura and Sheena (2005).) Therefore usually a large sample
approximation is employed (e.g.\ Theorem 13.5.1. of Anderson (2003)):
$$
\sqrt{n}\Bigl(\frac{l_i}{n}-\lambda_i\Bigr)\dlim N(0,2\lambda_i^2),\quad i=1,\ldots,p.
$$
Let $z_{\gamma}$ denote the upper $100\gamma$ percentile of the standard normal distribution. Since
\begin{align*}
&P\Bigl(\sqrt{\frac{n}{2}}\Bigl(\frac{l_1}{n\lambda_1}-1\Bigr)\geq z_{\gamma}\Bigr)\\
&=P\Bigl(l_1\geq (\sqrt{2n}z_{\gamma}+n)\lambda_1\Bigr)\rightarrow \gamma\quad \mbox{ as }n\rightarrow \infty,
\end{align*}
we have an approximate confidence interval
\begin{equation}
\label{old_interval}
\lambda_1\leq (\sqrt{2n}z_{\gamma}+n)^{-1}l_1,
\end{equation}
with confidence level close to $\gamma$ for sufficiently large $n$.

Now we propose an alternative approximation. Suppose $m=1$ in Theorem \ref{basic_theorem}, then as $\beta/\alpha$ goes to zero, 
$$
d_1\dlim \tilde{d}_1=\widetilde{\ww}_{11}.
$$
Since $\widetilde{\ww}_{11}/\xi_{1}\sim \chi^2(n)$,
$$
\frac{l_1}{\lambda_1}=\frac{d_1}{\xi_1}\dlim\chi^2(n).
$$
as $\beta/\alpha$ goes to zero.
{}From this asymptotics, we can make an approximate interval 
\begin{equation}
\label{new_interval}
\lambda_1\leq (\chi^2_{\gamma}(n))^{-1}l_1,
\end{equation}
where $\chi^2_{\gamma}(n)$ is the upper 100$\gamma$ percentile of $\chi^2$ distribution with the degree of freedom $n$. The interval (\ref{new_interval}) has approximately $\gamma$ confidence level when $\beta/\alpha$ is sufficiently close to zero.

We are interested in how large $n$ for  (\ref{old_interval}) or how
small $\beta/\alpha$ for (\ref{new_interval}) is required to get
practically sufficient approximations. Because of difficulty in
theoretical evaluations, we carried out a simulation study with the
fixed parameters $p=3,\ m=1,\ \xi_1=\xi_2=\xi_3=1,\ \alpha=1$, while
we select different $n$'s (5, 10, 20, 50, 100, 500, 1000) and
$\beta$'s (1.0, 0.9, 0.8, 0.6, 0.5, 0.3, 0.1, 0.01, 0.001). For each
case, 50000 Wishart random matrices are generated. We present the
results in Table 1. The numbers under L1(U1) indicate the ratio of the
largest sample eigenvalue which fall within the interval (\ref{old_interval})
with $\gamma=0.95(0.05)$, while those under U2(L2) show the similar
ratio with respect to the interval (\ref{new_interval}). Numbers in
bold indicate that they are within $\pm 0.01$ deviation from the desired value, hence the approximation may be good enough for many practical purposes.  We can summarize the result as follows;
\begin{table}[tbp]
\caption{Approximated Interval Estimation} \label{Table 1}
\begin{center}
\small
\begin{tabular}{lrrrrrrrrrrrr} \hline
\multicolumn{1}{|l|}{$\beta$} & 1\phantom{.}\phantom{000} &  &  & \multicolumn{1}{l|}{} & 0.9\phantom{00} &  &  & \multicolumn{1}{l|}{} & 0.8\phantom{00} &  &  & \multicolumn{1}{l|}{} \\ \hline
\multicolumn{1}{|l|}{} & \multicolumn{1}{l|}{U1} & \multicolumn{1}{l|}{U2} & \multicolumn{1}{l|}{L1} & \multicolumn{1}{l|}{L2} & \multicolumn{1}{l|}{U1} & \multicolumn{1}{l|}{U2} & \multicolumn{1}{l|}{L1} & \multicolumn{1}{l|}{L2} & \multicolumn{1}{l|}{U1} & \multicolumn{1}{l|}{U2} & \multicolumn{1}{l|}{L1} & \multicolumn{1}{l|}{L2} \\ \hline
\multicolumn{1}{|l|}{$n=5$} & \multicolumn{1}{r|}{.402} & \multicolumn{1}{r|}{.322} & \multicolumn{1}{r|}{1.00} & \multicolumn{1}{r|}{1.00} & \multicolumn{1}{r|}{.341} & \multicolumn{1}{r|}{.267} & \multicolumn{1}{r|}{1.00} & \multicolumn{1}{r|}{1.00} & \multicolumn{1}{r|}{.276} & \multicolumn{1}{r|}{.209} & \multicolumn{1}{r|}{1.00} & \multicolumn{1}{r|}{1.00} \\ \hline
\multicolumn{1}{|l|}{$n=10$} & \multicolumn{1}{r|}{.416} & \multicolumn{1}{r|}{.345} & \multicolumn{1}{r|}{1.00} & \multicolumn{1}{r|}{1.00} & \multicolumn{1}{r|}{.327} & \multicolumn{1}{r|}{.264} & \multicolumn{1}{r|}{1.00} & \multicolumn{1}{r|}{1.00} & \multicolumn{1}{r|}{.247} & \multicolumn{1}{r|}{.192} & \multicolumn{1}{r|}{1.00} & \multicolumn{1}{r|}{1.00} \\ \hline
\multicolumn{1}{|l|}{$n=20$} & \multicolumn{1}{r|}{.419} & \multicolumn{1}{r|}{.361} & \multicolumn{1}{r|}{1.00} & \multicolumn{1}{r|}{1.00} & \multicolumn{1}{r|}{.301} & \multicolumn{1}{r|}{.250} & \multicolumn{1}{r|}{1.00} & \multicolumn{1}{r|}{1.00} & \multicolumn{1}{r|}{.207} & \multicolumn{1}{r|}{.165} & \multicolumn{1}{r|}{1.00} & \multicolumn{1}{r|}{1.00} \\ \hline
\multicolumn{1}{|l|}{$n=50$} & \multicolumn{1}{r|}{.416} & \multicolumn{1}{r|}{.373} & \multicolumn{1}{r|}{1.00} & \multicolumn{1}{r|}{1.00} & \multicolumn{1}{r|}{.256} & \multicolumn{1}{r|}{.222} & \multicolumn{1}{r|}{1.00} & \multicolumn{1}{r|}{1.00} & \multicolumn{1}{r|}{.153} & \multicolumn{1}{r|}{.129} & \multicolumn{1}{r|}{1.00} & \multicolumn{1}{r|}{1.00} \\ \hline
\multicolumn{1}{|l|}{$n=100$} & \multicolumn{1}{r|}{.413} & \multicolumn{1}{r|}{.380} & \multicolumn{1}{r|}{1.00} & \multicolumn{1}{r|}{1.00} & \multicolumn{1}{r|}{.212} & \multicolumn{1}{r|}{.189} & \multicolumn{1}{r|}{1.00} & \multicolumn{1}{r|}{1.00} & \multicolumn{1}{r|}{.124} & \multicolumn{1}{r|}{.109} & \multicolumn{1}{r|}{.999} & \multicolumn{1}{r|}{.999} \\ \hline
\multicolumn{1}{|l|}{$n=500$} & \multicolumn{1}{r|}{.405} & \multicolumn{1}{r|}{.389} & \multicolumn{1}{r|}{1.00} & \multicolumn{1}{r|}{1.00} & \multicolumn{1}{r|}{.118} & \multicolumn{1}{r|}{.111} & \multicolumn{1}{r|}{.999} & \multicolumn{1}{r|}{.999} & \multicolumn{1}{r|}{.080} & \multicolumn{1}{r|}{.074} & \multicolumn{1}{r|}{.984} & \multicolumn{1}{r|}{.981} \\ \hline
\multicolumn{1}{|l|}{$n=1000$} & \multicolumn{1}{r|}{.407} & \multicolumn{1}{r|}{.395} & \multicolumn{1}{r|}{1.00} & \multicolumn{1}{r|}{1.00} & \multicolumn{1}{r|}{.097} & \multicolumn{1}{r|}{.093} & \multicolumn{1}{r|}{.995} & \multicolumn{1}{r|}{.994} & \multicolumn{1}{r|}{.070} & \multicolumn{1}{r|}{.067} & \multicolumn{1}{r|}{.972} & \multicolumn{1}{r|}{.970} \\ \hline
 &  &  &  &  &  &  &  &  &  &  &  &  \\ \hline
\multicolumn{1}{|l|}{$\beta$} & .6\phantom{00} &  &  & \multicolumn{1}{l|}{} & .5\phantom{00} &  &  & \multicolumn{1}{l|}{} & .3\phantom{00} &  &  & \multicolumn{1}{l|}{} \\ \hline
\multicolumn{1}{|l|}{} & \multicolumn{1}{l|}{U1} & \multicolumn{1}{l|}{U2} & \multicolumn{1}{l|}{L1} & \multicolumn{1}{l|}{L2} & \multicolumn{1}{l|}{U1} & \multicolumn{1}{l|}{U2} & \multicolumn{1}{l|}{L1} & \multicolumn{1}{l|}{L2} & \multicolumn{1}{l|}{U1} & \multicolumn{1}{l|}{U2} & \multicolumn{1}{l|}{L1} & \multicolumn{1}{l|}{L2} \\ \hline
\multicolumn{1}{|l|}{$n=5$} & \multicolumn{1}{r|}{.167} & \multicolumn{1}{r|}{.120} & \multicolumn{1}{r|}{1.00} & \multicolumn{1}{r|}{1.00} & \multicolumn{1}{r|}{.133} & \multicolumn{1}{r|}{.094} & \multicolumn{1}{r|}{1.00} & \multicolumn{1}{r|}{1.00} & \multicolumn{1}{r|}{.095} & \multicolumn{1}{r|}{.067} & \multicolumn{1}{r|}{1.00} & \multicolumn{1}{r|}{.998} \\ \hline
\multicolumn{1}{|l|}{$n=10$} & \multicolumn{1}{r|}{.139} & \multicolumn{1}{r|}{.104} & \multicolumn{1}{r|}{1.00} & \multicolumn{1}{r|}{1.00} & \multicolumn{1}{r|}{.111} & \multicolumn{1}{r|}{.083} & \multicolumn{1}{r|}{1.00} & \multicolumn{1}{r|}{.999} & \multicolumn{1}{r|}{.084} & \multicolumn{1}{r|}{.063} & \multicolumn{1}{r|}{.999} & \multicolumn{1}{r|}{.989} \\ \hline
\multicolumn{1}{|l|}{$n=20$} & \multicolumn{1}{r|}{.112} & \multicolumn{1}{r|}{.089} & \multicolumn{1}{r|}{1.00} & \multicolumn{1}{r|}{.998} & \multicolumn{1}{r|}{.094} & \multicolumn{1}{r|}{.073} & \multicolumn{1}{r|}{.999} & \multicolumn{1}{r|}{.995} & \multicolumn{1}{r|}{.075} & \multicolumn{1}{r|}{\bf.058} & \multicolumn{1}{r|}{.990} & \multicolumn{1}{r|}{.974} \\ \hline
\multicolumn{1}{|l|}{$n=50$} & \multicolumn{1}{r|}{.089} & \multicolumn{1}{r|}{.075} & \multicolumn{1}{r|}{.995} & \multicolumn{1}{r|}{.991} & \multicolumn{1}{r|}{.080} & \multicolumn{1}{r|}{.068} & \multicolumn{1}{r|}{.988} & \multicolumn{1}{r|}{.980} & \multicolumn{1}{r|}{.069} & \multicolumn{1}{r|}{\bf.058} & \multicolumn{1}{r|}{.974} & \multicolumn{1}{r|}{.961} \\ \hline
\multicolumn{1}{|l|}{$n=100$} & \multicolumn{1}{r|}{.078} & \multicolumn{1}{r|}{.068} & \multicolumn{1}{r|}{.984} & \multicolumn{1}{r|}{.978} & \multicolumn{1}{r|}{.072} & \multicolumn{1}{r|}{.064} & \multicolumn{1}{r|}{.975} & \multicolumn{1}{r|}{.968} & \multicolumn{1}{r|}{.062} & \multicolumn{1}{r|}{\bf.055} & \multicolumn{1}{r|}{.965} & \multicolumn{1}{r|}{\bf.957} \\ \hline
\multicolumn{1}{|l|}{$n=500$} & \multicolumn{1}{r|}{.064} & \multicolumn{1}{r|}{\bf.059} & \multicolumn{1}{r|}{.963} & \multicolumn{1}{r|}{\bf.959} & \multicolumn{1}{r|}{\bf.059} & \multicolumn{1}{r|}{\bf.055} & \multicolumn{1}{r|}{\bf.959} & \multicolumn{1}{r|}{\bf.955} & \multicolumn{1}{r|}{\bf.057} & \multicolumn{1}{r|}{\bf.053} & \multicolumn{1}{r|}{\bf.958} & \multicolumn{1}{r|}{\bf.954} \\ \hline
\multicolumn{1}{|l|}{$n=1000$} & \multicolumn{1}{r|}{\bf.060} & \multicolumn{1}{r|}{\bf.057} & \multicolumn{1}{r|}{\bf.960} & \multicolumn{1}{r|}{\bf.958} & \multicolumn{1}{r|}{\bf.058} & \multicolumn{1}{r|}{\bf.055} & \multicolumn{1}{r|}{\bf.956} & \multicolumn{1}{r|}{\bf.953} & \multicolumn{1}{r|}{\bf.055} & \multicolumn{1}{r|}{\bf.053} & \multicolumn{1}{r|}{\bf.956} & \multicolumn{1}{r|}{\bf.954} \\ \hline
 &  &  &  &  &  &  &  &  &  &  &  &  \\ \hline
\multicolumn{1}{|l|}{$\beta$} & .1\phantom{00} &  &  & \multicolumn{1}{l|}{} & .01\phantom{0} &  &  & \multicolumn{1}{l|}{} & .001 &  &  & \multicolumn{1}{l|}{} \\ \hline
\multicolumn{1}{|l|}{} & \multicolumn{1}{l|}{U1} & \multicolumn{1}{l|}{U2} & \multicolumn{1}{l|}{L1} & \multicolumn{1}{l|}{L2} & \multicolumn{1}{l|}{U1} & \multicolumn{1}{l|}{U2} & \multicolumn{1}{l|}{L1} & \multicolumn{1}{l|}{L2} & \multicolumn{1}{l|}{U1} & \multicolumn{1}{l|}{U2} & \multicolumn{1}{l|}{L1} & \multicolumn{1}{l|}{L2} \\ \hline
\multicolumn{1}{|l|}{$n=5$} & \multicolumn{1}{r|}{.075} & \multicolumn{1}{r|}{\bf.055} & \multicolumn{1}{r|}{1.00} & \multicolumn{1}{r|}{.975} & \multicolumn{1}{r|}{.068} & \multicolumn{1}{r|}{\bf.049} & \multicolumn{1}{r|}{1.00} & \multicolumn{1}{r|}{\bf.953} & \multicolumn{1}{r|}{.071} & \multicolumn{1}{r|}{\bf.052} & \multicolumn{1}{r|}{1.00} & \multicolumn{1}{r|}{\bf.951} \\ \hline
\multicolumn{1}{|l|}{$n=10$} & \multicolumn{1}{r|}{.072} & \multicolumn{1}{r|}{\bf.055} & \multicolumn{1}{r|}{.992} & \multicolumn{1}{r|}{\bf.959} & \multicolumn{1}{r|}{.068} & \multicolumn{1}{r|}{\bf.052} & \multicolumn{1}{r|}{.989} & \multicolumn{1}{r|}{\bf.951} & \multicolumn{1}{r|}{.065} & \multicolumn{1}{r|}{\bf.048} & \multicolumn{1}{r|}{.988} & \multicolumn{1}{r|}{\bf.950} \\ \hline
\multicolumn{1}{|l|}{$n=20$} & \multicolumn{1}{r|}{.066} & \multicolumn{1}{r|}{\bf.053} & \multicolumn{1}{r|}{.979} & \multicolumn{1}{r|}{\bf.956} & \multicolumn{1}{r|}{.066} & \multicolumn{1}{r|}{\bf.052} & \multicolumn{1}{r|}{.976} & \multicolumn{1}{r|}{\bf.951} & \multicolumn{1}{r|}{.063} & \multicolumn{1}{r|}{\bf.050} & \multicolumn{1}{r|}{.974} & \multicolumn{1}{r|}{\bf.948} \\ \hline
\multicolumn{1}{|l|}{$n=50$} & \multicolumn{1}{r|}{.061} & \multicolumn{1}{r|}{\bf.051} & \multicolumn{1}{r|}{.966} & \multicolumn{1}{r|}{\bf.953} & \multicolumn{1}{r|}{\bf.058} & \multicolumn{1}{r|}{\bf.048} & \multicolumn{1}{r|}{.964} & \multicolumn{1}{r|}{\bf.950} & \multicolumn{1}{r|}{\bf.059} & \multicolumn{1}{r|}{\bf.050} & \multicolumn{1}{r|}{.965} & \multicolumn{1}{r|}{\bf.950} \\ \hline
\multicolumn{1}{|l|}{$n=100$} & \multicolumn{1}{r|}{\bf.057} & \multicolumn{1}{r|}{\bf.050} & \multicolumn{1}{r|}{.961} & \multicolumn{1}{r|}{\bf.951} & \multicolumn{1}{r|}{\bf.057} & \multicolumn{1}{r|}{\bf.051} & \multicolumn{1}{r|}{.961} & \multicolumn{1}{r|}{\bf.951} & \multicolumn{1}{r|}{\bf.056} & \multicolumn{1}{r|}{\bf.049} & \multicolumn{1}{r|}{\bf.960} & \multicolumn{1}{r|}{\bf.951} \\ \hline
\multicolumn{1}{|l|}{$n=500$} & \multicolumn{1}{r|}{\bf.055} & \multicolumn{1}{r|}{\bf.051} & \multicolumn{1}{r|}{\bf.954} & \multicolumn{1}{r|}{\bf.950} & \multicolumn{1}{r|}{\bf.054} & \multicolumn{1}{r|}{\bf.050} & \multicolumn{1}{r|}{\bf.956} & \multicolumn{1}{r|}{\bf.952} & \multicolumn{1}{r|}{\bf.053} & \multicolumn{1}{r|}{\bf.049} & \multicolumn{1}{r|}{\bf.954} & \multicolumn{1}{r|}{\bf.950} \\ \hline
\multicolumn{1}{|l|}{$n=1000$} & \multicolumn{1}{r|}{\bf.052} & \multicolumn{1}{r|}{\bf.050} & \multicolumn{1}{r|}{\bf.954} & \multicolumn{1}{r|}{\bf.951} & \multicolumn{1}{r|}{\bf.053} & \multicolumn{1}{r|}{\bf.050} & \multicolumn{1}{r|}{\bf.953} & \multicolumn{1}{r|}{\bf.949} & \multicolumn{1}{r|}{\bf.052} & \multicolumn{1}{r|}{\bf.049} & \multicolumn{1}{r|}{\bf.954} & \multicolumn{1}{r|}{\bf.951} \\ \hline
\end{tabular}
\end{center}
\end{table} %
\begin{enumerate}
\item In every case (\ref{new_interval}) gives better approximation than (\ref{old_interval}).
\item Since $\beta$ is as large as 0.3, (\ref{new_interval}) already gives a good approximation. In that sense, the approximated interval (\ref{new_interval}) seems robust. When $\beta$ is smaller or equal to 0.1, (\ref{new_interval}) works well even with small samples such as $n=5$ or 10, while (\ref{old_interval}) needs samples as large as 100 or 500. 
\item  When $\beta$ is from 0.5 to 0.6,  both approximations need large samples such as 500 or 1000. If $\beta$ is larger than 0.6, they need samples larger than 1000 for a good approximation.
\end{enumerate}
Similarly we can make an interval estimation for the smallest eigenvalue, $\lambda_p$. Let $m=p-1$ in Theorem \ref{basic_theorem}, then  $\beta/\alpha$ goes to zero, 
$$
d_p\dlim \tilde{d}_p=\widetilde{\ww}_{22}.
$$
Since $\widetilde{\ww}_{22}/\xi_{p}\sim \chi^2(n-p+1)$,
$$
\frac{l_p}{\lambda_p}=\frac{d_p}{\xi_p}\dlim\chi^2(n-p+1).
$$
as $\beta/\alpha$ goes to zero. Using this fact, we can estimate $\lambda_p$ to lie in the interval
\begin{equation}
\label{new_interval_2}
\lambda_p\leq (\chi^2_{\gamma}(n-p+1))^{-1}l_p,
\end{equation}
at approximately $\gamma$ confidence level when $\beta/\alpha$ is sufficiently close to zero.

We now compare (\ref{old_interval}) and  (\ref{new_interval})
more closely in view of the known results on asymptotic expansion of
distribution of sample eigenvalues.
Let
\begin{equation}
\label{def_A_n}
A_n=\sqrt{\frac{n}{2}}\left(\frac{l_1/\lambda_1}{n}-1\right),
\end{equation}
The asymptotic expansion of $A_n$ up to the order $n^{1/2}$ is given by
(see Sugiura (1973))
\begin{equation}
\label{asym_dist_A_n}
F_{A_n}(t)=\Phi(t)-\frac{\sqrt{2}\phi(t)}{3\sqrt{n}}\left((t^2-1)+\frac{1}{2}\sum_{i=2}^p
\frac{\lambda_i}{\lambda_1-\lambda_i}\right)+o(n^{-1/2}).
\end{equation}
Now suppose $x_i,\ i=1,\ldots,n$, are independently  and identically distributed as $\chi^2(1)$ distribution.
The normalized variable 
$$
\tilde{x}_i=\frac{1}{\sqrt{2}}(x_i-1),\quad i=1,\ldots,n
$$ 
has zero mean and unit variance. Let
\begin{equation}
\label{def_B_n}
B_n=\frac{1}{\sqrt{n}}\sum_{i=1}^n\tilde{x}_i=\sqrt{\frac{n}{2}}\left(\frac{\sum_{i=1}^nx_i}{n}-1\right).
\end{equation}
The asymptotic expansion of $B_n$ up to the order $n^{1/2}$ is given by
\begin{equation}
\label{asym_dist_B_n}
F_{B_n}(t)=\Phi(t)-\frac{\sqrt{2}\phi(t)}{3\sqrt{n}}(t^2-1)+o(n^{-1/2}).
\end{equation}

Comparing (\ref{asym_dist_A_n}) and (\ref{asym_dist_B_n}), we notice that if $t^2>1$, then the absolute value of the second term in (\ref{asym_dist_B_n}) is smaller than that of (\ref{asym_dist_A_n}) by the margin
\begin{equation}
\label{asym_margin}
\frac{1}{2}\sum_{i=2}^p\frac{\lambda_i}{\lambda_1-\lambda_i}
\end{equation}
Since $l_1/\lambda_1$ is asymptotically distributed as $\chi^2(n)$
when the largest population eigenvalue is infinitely deviated from the
others, $l_1/\lambda_1$ in (\ref{def_A_n}) is similarly distributed as
$\sum x_i$ in (\ref{def_B_n}). In this case (\ref{asym_margin})
vanishes and both expansions (\ref{asym_dist_A_n}),
(\ref{asym_dist_B_n}) coincide. It is naturally conjectured that when
the largest population eigenvalue $\lambda_1$ is positioned far away
from the smaller eigenvalues, we can make an ``easier'' inference on
$\lambda_1$. The fact that the term (\ref{asym_margin}) shrinks in
that situation supports this conjecture as well as our simulation results.

%
%
%
%
%
%
%
%
%
\subsection{Testing Equality of the Smallest Eigenvalues}
As in the introduction of Section 11.7.3 of Anderson (2003), the equality of the $p-m$ smallest population eigenvalues
\begin{equation}
\label{equ_small_eigen}
\lambda_{m+1}=\cdots=\lambda_p\quad(\mbox{say }\sigma^2)
\end{equation}
is equivalent to the covariance structure 
$$
\s=\bm{\Phi}+\sigma^2\iip,
$$
where $\bm{\Phi}$, a positive semidefinite matrix with rank $m$,
represents the variance-covariance matrix of a systematic part and
$\sigma^2\iip$ arises from measurement error. If hypothesis
(\ref{equ_small_eigen}) is accepted, then it suggests that the
systematic part might consist of $m$ independent factors. Need for
testing (\ref{equ_small_eigen}) also arises in principal component
analysis when the dimension of principal components has to be decided.
Once it is accepted and $\sigma^2$ is sufficiently small, which might
require another hypothesis testing, we could ignore the last $p-m$
principal components.

The likelihood ratio statistic for testing (\ref{equ_small_eigen}) is
given (see e.g.\ Theorem 9.6.1 of Muirhead (1982)) by
$$
\bm{V}=\frac{\prod_{i=m+1}^p l_i}{\Bigl(\sum_{i=m+1}^p l_i\Bigr)^{p-m}}(p-m)^{p-m},
$$
and the critical region is $\bm{V}\leq c(\gamma)$ for a given significance level $\gamma$ . 

In order to give the critical point $c(\gamma)$, we traditionally make use of the asymptotic convergence
\begin{equation}
\label{eq:likelihood-ratio-V}
-n\log \bm{V}\dlim \chi^2((p-m+2)(p-m-1)/2),\quad \mbox{as $n\rightarrow \infty$.}
\end{equation}
Bartlett adjustment and further refinement on the asymptotic result
are found in Section 9.6 of Muirhead (1982).  {}From this convergence,
the approximate critical point is given as
\begin{equation}
\label{criti_point_1}
c(\gamma)=\exp\Bigl(-n^{-1}\chi^2_{\gamma}((p-m+2)(p-m-1)/2)\Bigr)
\end{equation}

On the other hand we can approximate the critical point $c(\gamma)$
based on the asymptotic result in Theorem \ref{basic_theorem}.  We can
expect that this approach yields good approximation since in testing
hypothesis (\ref{equ_small_eigen}), we often encounter the situation
where the eigenvalues $\lambda_{m+1},\ldots,\lambda_p$ are much
smaller than the other eigenvalues. The hypothesis
(\ref{equ_small_eigen}) with small $\sigma^2$ corresponds to the case
$\xi_{m+1}=\cdots=\xi_{p}=1$ in Theorem
\ref{basic_theorem}. Consequently we can approximate the distribution
of $V$ in (\ref{eq:likelihood-ratio-V}) by the distribution of
$$
\widetilde{\bm{V}}=\frac{\prod_{i=m+1}^p d_i}{\Bigl(\sum_{i=m+1}^p d_i\Bigr)^{p-m}}(p-m)^{p-m},
$$
where $d_i\  (i=m+1,\ldots,p)$ are the eigenvalues of Wishart matrix
$\ww_{p-m}(n-m,\bm{I}_{p-m}).$ Even under the distribution
$\ww_{p-m}(n-m,\bm{I}_{p-m})$, it is not easy to derive analytical
expressions for percentage points for $\widetilde{\bm{V}}$. For $p-m=2$,  the distribution function  is explicitly given (see 10.7.3. of Anderson (2003)) by
$$
P(\widetilde{\bm{V}}\leq v)=v^{(n-m-1)/2}
$$
which gives the critical point $c(\gamma)$ as
\begin{equation}
\label{criti_point_2}
c(\gamma)=\gamma^{2/(n-m-1)}.
\end{equation}
Generally a numerical calculation is needed for the exact evaluation of critical points. For this problem, refer to Consul (1967) and Pillai and Nagarsenker (1971). 

We made a simulation for the comparison between the above two methods. Let $p=3,\ m=1$ and consider testing the hypothesis $\lambda_2=\lambda_3$. We examined the accuracy of the two critical points (\ref{criti_point_1}) and (\ref{criti_point_2}) 
with $\gamma=0.05$ and $\gamma=0.01$ by simulating the probability of
$\bm{V}$ being smaller than these critical points when
$\lambda_2=\lambda_3$, that is, the probability of the error of the
first kind. We put $\lambda_1=1,\ \lambda_2=\beta,\ \lambda_3=\beta$
and varied both $\beta$ and $n.$ Table 2 shows the result, where the
labels 5(1)\%1 and 5(1)\%2 indicate that the numbers below correspond
to the critical points (\ref{criti_point_1}) and (\ref{criti_point_2})
respectively with $\gamma=0.05(0.01).$  Numbers in bold mean that
they are within $\pm0.001$ deviation from the desired value.

\begin{table}[tbp]
\caption{Simulated Type 1 Error}\label{Table 2}
\begin{center}
\small
\begin{tabular}{lrrrrrrrrrrrr} \hline
\multicolumn{1}{|l|}{$\beta$} & 1\phantom{.}\phantom{000} &  &  & \multicolumn{1}{l|}{} & .9\phantom{00} &  &  & \multicolumn{1}{l|}{} & .8\phantom{00} &  &  & \multicolumn{1}{l|}{} \\ \hline
\multicolumn{1}{|l|}{} & \multicolumn{1}{l|}{5\%1} & \multicolumn{1}{l|}{5\%2} & \multicolumn{1}{l|}{1\%1} & \multicolumn{1}{l|}{1\%2} & \multicolumn{1}{l|}{5\%1} & \multicolumn{1}{l|}{5\%2} & \multicolumn{1}{l|}{1\%1} & \multicolumn{1}{l|}{1\%2} & \multicolumn{1}{l|}{5\%1} & \multicolumn{1}{l|}{5\%2} & \multicolumn{1}{l|}{1\%1} & \multicolumn{1}{l|}{1\%2} \\ \hline
\multicolumn{1}{|l|}{$n=5$} & \multicolumn{1}{r|}{.140} & \multicolumn{1}{r|}{.041} & \multicolumn{1}{r|}{.051} & \multicolumn{1}{r|}{.008} & \multicolumn{1}{r|}{.142} & \multicolumn{1}{r|}{.041} & \multicolumn{1}{r|}{.053} & \multicolumn{1}{r|}{.008} & \multicolumn{1}{r|}{.141} & \multicolumn{1}{r|}{.042} & \multicolumn{1}{r|}{.054} & \multicolumn{1}{r|}{.008} \\ \hline
\multicolumn{1}{|l|}{$n=10$} & \multicolumn{1}{r|}{.063} & \multicolumn{1}{r|}{.033} & \multicolumn{1}{r|}{.015} & \multicolumn{1}{r|}{.005} & \multicolumn{1}{r|}{.064} & \multicolumn{1}{r|}{.033} & \multicolumn{1}{r|}{.016} & \multicolumn{1}{r|}{.006} & \multicolumn{1}{r|}{.067} & \multicolumn{1}{r|}{.036} & \multicolumn{1}{r|}{.017} & \multicolumn{1}{r|}{.007} \\ \hline
\multicolumn{1}{|l|}{$n=20$} & \multicolumn{1}{r|}{.038} & \multicolumn{1}{r|}{.026} & \multicolumn{1}{r|}{.007} & \multicolumn{1}{r|}{.004} & \multicolumn{1}{r|}{.039} & \multicolumn{1}{r|}{.028} & \multicolumn{1}{r|}{.008} & \multicolumn{1}{r|}{.005} & \multicolumn{1}{r|}{.041} & \multicolumn{1}{r|}{.030} & \multicolumn{1}{r|}{.008} & \multicolumn{1}{r|}{.005} \\ \hline
\multicolumn{1}{|l|}{$n=50$} & \multicolumn{1}{r|}{.025} & \multicolumn{1}{r|}{.021} & \multicolumn{1}{r|}{.004} & \multicolumn{1}{r|}{.003} & \multicolumn{1}{r|}{.027} & \multicolumn{1}{r|}{.024} & \multicolumn{1}{r|}{.004} & \multicolumn{1}{r|}{.003} & \multicolumn{1}{r|}{.032} & \multicolumn{1}{r|}{.029} & \multicolumn{1}{r|}{.006} & \multicolumn{1}{r|}{.005} \\ \hline
\multicolumn{1}{|l|}{$n=100$} & \multicolumn{1}{r|}{.022} & \multicolumn{1}{r|}{.020} & \multicolumn{1}{r|}{.003} & \multicolumn{1}{r|}{.003} & \multicolumn{1}{r|}{.025} & \multicolumn{1}{r|}{.023} & \multicolumn{1}{r|}{.003} & \multicolumn{1}{r|}{.003} & \multicolumn{1}{r|}{.034} & \multicolumn{1}{r|}{.031} & \multicolumn{1}{r|}{.005} & \multicolumn{1}{r|}{.005} \\ \hline
\multicolumn{1}{|l|}{$n=500$} & \multicolumn{1}{r|}{.016} & \multicolumn{1}{r|}{.016} & \multicolumn{1}{r|}{.002} & \multicolumn{1}{r|}{.002} & \multicolumn{1}{r|}{.029} & \multicolumn{1}{r|}{.029} & \multicolumn{1}{r|}{.005} & \multicolumn{1}{r|}{.005} & \multicolumn{1}{r|}{.043} & \multicolumn{1}{r|}{.043} & \multicolumn{1}{r|}{.008} & \multicolumn{1}{r|}{.008} \\ \hline
\multicolumn{1}{|l|}{$n=1000$} & \multicolumn{1}{r|}{.017} & \multicolumn{1}{r|}{.017} & \multicolumn{1}{r|}{.002} & \multicolumn{1}{r|}{.002} & \multicolumn{1}{r|}{.035} & \multicolumn{1}{r|}{.035} & \multicolumn{1}{r|}{.006} & \multicolumn{1}{r|}{.006} & \multicolumn{1}{r|}{.048} & \multicolumn{1}{r|}{.047} & \multicolumn{1}{r|}{\bf.009} & \multicolumn{1}{r|}{\bf.009} \\ \hline
 &  &  &  &  &  &  &  &  &  &  &  &  \\ \hline
\multicolumn{1}{|l|}{$\beta$} & .6\phantom{00} &  &  & \multicolumn{1}{l|}{} & .5\phantom{00} &  &  & \multicolumn{1}{l|}{} & .3\phantom{00} &  &  & \multicolumn{1}{l|}{} \\ \hline
\multicolumn{1}{|l|}{} & \multicolumn{1}{l|}{5\%1} & \multicolumn{1}{l|}{5\%2} & \multicolumn{1}{l|}{1\%1} & \multicolumn{1}{l|}{1\%2} & \multicolumn{1}{l|}{5\%1} & \multicolumn{1}{l|}{5\%2} & \multicolumn{1}{l|}{1\%1} & \multicolumn{1}{l|}{1\%2} & \multicolumn{1}{l|}{5\%1} & \multicolumn{1}{l|}{5\%2} & \multicolumn{1}{l|}{1\%1} & \multicolumn{1}{l|}{1\%2} \\ \hline
\multicolumn{1}{|l|}{$n=5$} & \multicolumn{1}{r|}{.145} & \multicolumn{1}{r|}{.044} & \multicolumn{1}{r|}{.055} & \multicolumn{1}{r|}{.008} & \multicolumn{1}{r|}{.148} & \multicolumn{1}{r|}{.044} & \multicolumn{1}{r|}{.056} & \multicolumn{1}{r|}{\bf.009} & \multicolumn{1}{r|}{.158} & \multicolumn{1}{r|}{.047} & \multicolumn{1}{r|}{.059} & \multicolumn{1}{r|}{\bf.009} \\ \hline
\multicolumn{1}{|l|}{$n=10$} & \multicolumn{1}{r|}{.072} & \multicolumn{1}{r|}{.039} & \multicolumn{1}{r|}{.019} & \multicolumn{1}{r|}{.008} & \multicolumn{1}{r|}{.076} & \multicolumn{1}{r|}{.041} & \multicolumn{1}{r|}{.021} & \multicolumn{1}{r|}{.008} & \multicolumn{1}{r|}{.086} & \multicolumn{1}{r|}{.046} & \multicolumn{1}{r|}{.023} & \multicolumn{1}{r|}{\bf.009} \\ \hline
\multicolumn{1}{|l|}{$n=20$} & \multicolumn{1}{r|}{\bf.051} & \multicolumn{1}{r|}{.036} & \multicolumn{1}{r|}{\bf.011} & \multicolumn{1}{r|}{.006} & \multicolumn{1}{r|}{.057} & \multicolumn{1}{r|}{.042} & \multicolumn{1}{r|}{.012} & \multicolumn{1}{r|}{.005} & \multicolumn{1}{r|}{.067} & \multicolumn{1}{r|}{\bf.049} & \multicolumn{1}{r|}{.016} & \multicolumn{1}{r|}{\bf.009} \\ \hline
\multicolumn{1}{|l|}{$n=50$} & \multicolumn{1}{r|}{.046} & \multicolumn{1}{r|}{.040} & \multicolumn{1}{r|}{\bf.009} & \multicolumn{1}{r|}{.008} & \multicolumn{1}{r|}{.053} & \multicolumn{1}{r|}{.047} & \multicolumn{1}{r|}{\bf.011} & \multicolumn{1}{r|}{\bf.009} & \multicolumn{1}{r|}{.056} & \multicolumn{1}{r|}{\bf.050} & \multicolumn{1}{r|}{.012} & \multicolumn{1}{r|}{\bf.010} \\ \hline
\multicolumn{1}{|l|}{$n=100$} & \multicolumn{1}{r|}{\bf.050} & \multicolumn{1}{r|}{.047} & \multicolumn{1}{r|}{\bf.010} & \multicolumn{1}{r|}{\bf.009} & \multicolumn{1}{r|}{\bf.050} & \multicolumn{1}{r|}{.047} & \multicolumn{1}{r|}{\bf.010} & \multicolumn{1}{r|}{\bf.009} & \multicolumn{1}{r|}{.052} & \multicolumn{1}{r|}{\bf.049} & \multicolumn{1}{r|}{\bf.010} & \multicolumn{1}{r|}{\bf.009} \\ \hline
\multicolumn{1}{|l|}{$n=500$} & \multicolumn{1}{r|}{\bf.050} & \multicolumn{1}{r|}{\bf.050} & \multicolumn{1}{r|}{\bf.010} & \multicolumn{1}{r|}{\bf.009} & \multicolumn{1}{r|}{\bf.050} & \multicolumn{1}{r|}{\bf.050} & \multicolumn{1}{r|}{\bf.010} & \multicolumn{1}{r|}{\bf.010} & \multicolumn{1}{r|}{.052} & \multicolumn{1}{r|}{\bf.051} & \multicolumn{1}{r|}{\bf.010} & \multicolumn{1}{r|}{\bf.010} \\ \hline
\multicolumn{1}{|l|}{$n=1000$} & \multicolumn{1}{r|}{\bf.051} & \multicolumn{1}{r|}{\bf.050} & \multicolumn{1}{r|}{\bf.010} & \multicolumn{1}{r|}{\bf.010} & \multicolumn{1}{r|}{\bf.051} & \multicolumn{1}{r|}{\bf.050} & \multicolumn{1}{r|}{\bf.010} & \multicolumn{1}{r|}{\bf.010} & \multicolumn{1}{r|}{\bf.051} & \multicolumn{1}{r|}{\bf.051} & \multicolumn{1}{r|}{\bf.010} & \multicolumn{1}{r|}{\bf.010} \\ \hline
 &  &  &  &  &  &  &  &  &  &  &  &  \\ \hline
\multicolumn{1}{|l|}{$\beta$} & .1\phantom{00} &  &  & \multicolumn{1}{l|}{} & .01\phantom{0} &  &  & \multicolumn{1}{l|}{} & .001 &  &  & \multicolumn{1}{l|}{} \\ \hline
\multicolumn{1}{|l|}{} & \multicolumn{1}{l|}{5\%1} & \multicolumn{1}{l|}{5\%2} & \multicolumn{1}{l|}{1\%1} & \multicolumn{1}{l|}{1\%2} & \multicolumn{1}{l|}{5\%1} & \multicolumn{1}{l|}{5\%2} & \multicolumn{1}{l|}{1\%1} & \multicolumn{1}{l|}{1\%2} & \multicolumn{1}{l|}{5\%1} & \multicolumn{1}{l|}{5\%2} & \multicolumn{1}{l|}{1\%1} & \multicolumn{1}{l|}{1\%2} \\ \hline
\multicolumn{1}{|l|}{$n=5$} & \multicolumn{1}{r|}{.161} & \multicolumn{1}{r|}{\bf.049} & \multicolumn{1}{r|}{.061} & \multicolumn{1}{r|}{\bf.010} & \multicolumn{1}{r|}{.164} & \multicolumn{1}{r|}{\bf.050} & \multicolumn{1}{r|}{.063} & \multicolumn{1}{r|}{\bf.011} & \multicolumn{1}{r|}{.167} & \multicolumn{1}{r|}{\bf.051} & \multicolumn{1}{r|}{.064} & \multicolumn{1}{r|}{\bf.010} \\ \hline
\multicolumn{1}{|l|}{$n=10$} & \multicolumn{1}{r|}{.091} & \multicolumn{1}{r|}{\bf.051} & \multicolumn{1}{r|}{.026} & \multicolumn{1}{r|}{\bf.011} & \multicolumn{1}{r|}{.092} & \multicolumn{1}{r|}{\bf.051} & \multicolumn{1}{r|}{.025} & \multicolumn{1}{r|}{\bf.009} & \multicolumn{1}{r|}{.090} & \multicolumn{1}{r|}{\bf.050} & \multicolumn{1}{r|}{.024} & \multicolumn{1}{r|}{\bf.010} \\ \hline
\multicolumn{1}{|l|}{$n=20$} & \multicolumn{1}{r|}{.066} & \multicolumn{1}{r|}{\bf.049} & \multicolumn{1}{r|}{.015} & \multicolumn{1}{r|}{\bf.010} & \multicolumn{1}{r|}{.067} & \multicolumn{1}{r|}{\bf.050} & \multicolumn{1}{r|}{.016} & \multicolumn{1}{r|}{\bf.010} & \multicolumn{1}{r|}{.067} & \multicolumn{1}{r|}{\bf.050} & \multicolumn{1}{r|}{.016} & \multicolumn{1}{r|}{\bf.010} \\ \hline
\multicolumn{1}{|l|}{$n=50$} & \multicolumn{1}{r|}{.057} & \multicolumn{1}{r|}{\bf.051} & \multicolumn{1}{r|}{.012} & \multicolumn{1}{r|}{\bf.010} & \multicolumn{1}{r|}{.057} & \multicolumn{1}{r|}{\bf.051} & \multicolumn{1}{r|}{.012} & \multicolumn{1}{r|}{\bf.010} & \multicolumn{1}{r|}{.055} & \multicolumn{1}{r|}{.048} & \multicolumn{1}{r|}{.012} & \multicolumn{1}{r|}{\bf.010} \\ \hline
\multicolumn{1}{|l|}{$n=100$} & \multicolumn{1}{r|}{.053} & \multicolumn{1}{r|}{\bf.050} & \multicolumn{1}{r|}{\bf.011} & \multicolumn{1}{r|}{\bf.010} & \multicolumn{1}{r|}{.053} & \multicolumn{1}{r|}{\bf.049} & \multicolumn{1}{r|}{\bf.011} & \multicolumn{1}{r|}{\bf.010} & \multicolumn{1}{r|}{.054} & \multicolumn{1}{r|}{\bf.051} & \multicolumn{1}{r|}{\bf.011} & \multicolumn{1}{r|}{\bf.010} \\ \hline
\multicolumn{1}{|l|}{$n=500$} & \multicolumn{1}{r|}{\bf.051} & \multicolumn{1}{r|}{\bf.051} & \multicolumn{1}{r|}{\bf.010} & \multicolumn{1}{r|}{\bf.010} & \multicolumn{1}{r|}{\bf.050} & \multicolumn{1}{r|}{\bf.049} & \multicolumn{1}{r|}{\bf.010} & \multicolumn{1}{r|}{\bf.010} & \multicolumn{1}{r|}{\bf.051} & \multicolumn{1}{r|}{\bf.050} & \multicolumn{1}{r|}{\bf.011} & \multicolumn{1}{r|}{\bf.010} \\ \hline
\multicolumn{1}{|l|}{$n=1000$} & \multicolumn{1}{r|}{\bf.050} & \multicolumn{1}{r|}{\bf.050} & \multicolumn{1}{r|}{\bf.009} & \multicolumn{1}{r|}{\bf.009} & \multicolumn{1}{r|}{\bf.051} & \multicolumn{1}{r|}{\bf.051} & \multicolumn{1}{r|}{\bf.010} & \multicolumn{1}{r|}{\bf.010} & \multicolumn{1}{r|}{\bf.051} & \multicolumn{1}{r|}{\bf.050} & \multicolumn{1}{r|}{\bf.010} & \multicolumn{1}{r|}{\bf.010} \\ \hline
\end{tabular}
\end{center}
\end{table}
We can summarize the result as follows;
\begin{enumerate}
\item If $\beta\geq 0.8$, both (\ref{criti_point_1}) and
  (\ref{criti_point_2}) need a large sample size. Especially when $\beta$ is as large as 1.0 or 0.9, more than 1000 samples are required for a good approximation. There is no meaningful difference between both critical points.

\item If $\beta$ equals 0.6 or 0.5, 50 (sometimes 20) samples are large enough to give a good approximation for both (\ref{criti_point_1}) and (\ref{criti_point_2}).  There is no significant difference between both critical points.

\item If $\beta<0.5$, (\ref{criti_point_2}) shows significantly better performance than (\ref{criti_point_1}). Even with such a small sample as 5, (\ref{criti_point_2}) gives very accurate approximations. The critical point (\ref{criti_point_2}) is robust in the sense that it already gives an excellent approximation when the smallest eigenvalues are 0.3 times as large as the largest eigenvalue.   
\end{enumerate}
%
%
%
%
%
%
%
%

\end{document}